\documentclass{amsart}
\usepackage{amsfonts}

\newtheorem{theorem}{Theorem}
\theoremstyle{plain}

\newtheorem{lemma}{Lemma}

\newtheorem{remark}{Remark}

\numberwithin{equation}{section}
\input{tcilatex}

\begin{document}
\title[Upper Bounds for the Distance to Finite-Dimensional Subspaces]{Upper
Bounds for the Distance to Finite-Dimensional Subspaces in Inner Product
Spaces}
\author{S.S. Dragomir}
\address{School of Computer Science and Mathematics\\
Victoria University of Technology\\
PO Box 14428, MCMC 8001\\
VIC, Australia.}
\email{sever@matilda.vu.edu.au}
\urladdr{http://rgmia.vu.edu.au/SSDragomirWeb.html}
\date{22 October, 2004.}
\subjclass[2000]{46C05, 26D15.}
\keywords{Finite-dimensional subspaces, Distance, Bessel's inequality,
Boas-Bellman's inequality, Bombieri's inequality, Hadamard's inequality,
Gram's inequality.}

\begin{abstract}
We establish upper bounds for the distance to finite-dimensional subspaces
in inner product spaces and improve some generalisations of Bessel's
inequality obtained by Boas, Bellman and Bombieri. Refinements of the
Hadamard inequality for Gram determinants are also given.
\end{abstract}

\maketitle

\section{Introduction}

Let $\left( H;\left\langle \cdot ,\cdot \right\rangle \right) $ be an inner
product space over the real or complex number field $\mathbb{K}$, $\left\{
y_{1},\dots ,y_{n}\right\} $ a subset of $H$ and $G\left( y_{1},\dots
,y_{n}\right) $ the \textit{gram matrix} of $\left\{ y_{1},\dots
,y_{n}\right\} $ where $\left( i,j\right) -$entry is $\left\langle
y_{i},y_{j}\right\rangle .$ The determinant of $G\left( y_{1},\dots
,y_{n}\right) $ is called the \textit{Gram determinant} of $\left\{
y_{1},\dots ,y_{n}\right\} $ and is denoted by $\Gamma \left( y_{1},\dots
,y_{n}\right) .$ Thus,%
\begin{equation*}
\Gamma \left( y_{1},\dots ,y_{n}\right) =\left\vert 
\begin{array}{c}
\left\langle y_{1},y_{1}\right\rangle \;\left\langle
y_{1},y_{2}\right\rangle \;\cdots \;\left\langle y_{1},y_{n}\right\rangle \\ 
\left\langle y_{2},y_{1}\right\rangle \;\left\langle
y_{2},y_{2}\right\rangle \;\cdots \;\left\langle y_{2},y_{n}\right\rangle \\ 
\cdots \cdots \cdots \cdots \cdots \\ 
\left\langle y_{n},y_{1}\right\rangle \;\left\langle
y_{n},y_{2}\right\rangle \;\cdots \;\left\langle y_{n},y_{n}\right\rangle%
\end{array}%
\right\vert .
\end{equation*}

Following \cite[p. 129 -- 133]{DE}, we state here some general results for
the Gram determinant that will be used in the sequel.

\begin{enumerate}
\item Let $\left\{ x_{1},\dots ,x_{n}\right\} \subset H.$ Then $\Gamma
\left( x_{1},\dots ,x_{n}\right) \neq 0$ if and only if $\left\{ x_{1},\dots
,x_{n}\right\} $ is linearly independent;

\item Let $M=span\left\{ x_{1},\dots ,x_{n}\right\} $ be $n-$dimensional in $%
H,$ i.e., $\left\{ x_{1},\dots ,x_{n}\right\} $ is linearly independent.
Then for each $x\in H,$ the distance $d\left( x,M\right) $ from $x$ to the
linear subspace $H$ has the representations%
\begin{equation}
d^{2}\left( x,M\right) =\frac{\Gamma \left( x_{1},\dots ,x_{n},x\right) }{%
\Gamma \left( x_{1},\dots ,x_{n}\right) }  \label{1.1}
\end{equation}%
and%
\begin{equation}
d^{2}\left( x,M\right) =\left\Vert x\right\Vert ^{2}-\beta ^{T}G^{-1}\beta ,
\label{1.2}
\end{equation}%
where $G=G\left( x_{1},\dots ,x_{n}\right) ,$ $G^{-1}$ is the inverse matrix
of $G$ and 
\begin{equation*}
\beta ^{T}=\left( \left\langle x,x_{1}\right\rangle ,\left\langle
x,x_{2}\right\rangle ,\dots ,\left\langle x,x_{n}\right\rangle \right) ,
\end{equation*}%
denotes the transpose of the column vector $\beta .$

Moreover, one has the simpler representation%
\begin{equation}
d^{2}\left( x,M\right) =\left\{ 
\begin{array}{ll}
\left\Vert x\right\Vert ^{2}-\frac{\left( \sum_{i=1}^{n}\left\vert
\left\langle x,x_{i}\right\rangle \right\vert ^{2}\right) ^{2}}{\left\Vert
\sum_{i=1}^{n}\left\langle x,x_{i}\right\rangle x_{i}\right\Vert ^{2}} & 
\text{if \ }x\notin M^{\perp }, \\ 
&  \\ 
\left\Vert x\right\Vert ^{2} & \text{if \ }x\in M^{\perp },%
\end{array}%
\right.  \label{1.3}
\end{equation}%
where $M^{\perp }$ denotes the orthogonal complement of $M.$

\item Let $\left\{ x_{1},\dots ,x_{n}\right\} $ be a set of nonzero vectors
in $H.$ Then%
\begin{equation}
0\leq \Gamma \left( x_{1},\dots ,x_{n}\right) \leq \left\Vert
x_{1}\right\Vert ^{2}\left\Vert x_{2}\right\Vert ^{2}\cdots \left\Vert
x_{n}\right\Vert ^{2}.  \label{1.4}
\end{equation}%
The equality holds on the left (respectively right) side of (\ref{1.4}) if
and only if $\left\{ x_{1},\dots ,x_{n}\right\} $ is linearly dependent
(respectively orthogonal). The first inequality in (\ref{1.4}) is known in
the literature as \textit{Gram's inequality} while the second one is known
as \textit{Hadamard's inequality.}

\item If $\left\{ x_{1},\dots ,x_{n}\right\} $ is an orthonormal set in $H,$
i.e., $\left\langle x_{i},x_{j}\right\rangle =\delta _{ij},$ $i,j\in \left\{
1,\dots ,n\right\} ,$ where $\delta _{ij}$ is Kronecker's delta, then%
\begin{equation}
d^{2}\left( x,M\right) =\left\Vert x\right\Vert
^{2}-\sum_{i=1}^{n}\left\vert \left\langle x,x_{i}\right\rangle \right\vert
^{2}.  \label{1.5}
\end{equation}
\end{enumerate}

The following inequalities which involve Gram determinants may be stated as
well \cite[p. 597]{MPF}:%
\begin{equation}
\frac{\Gamma \left( x_{1},\dots ,x_{n}\right) }{\Gamma \left( x_{1},\dots
,x_{k}\right) }\leq \frac{\Gamma \left( x_{2},\dots ,x_{n}\right) }{\Gamma
\left( x_{1},\dots ,x_{k}\right) }\leq \cdots \leq \Gamma \left(
x_{k+1},\dots ,x_{n}\right) ,  \label{1.6}
\end{equation}%
\begin{equation}
\Gamma \left( x_{1},\dots ,x_{n}\right) \leq \Gamma \left( x_{1},\dots
,x_{k}\right) \Gamma \left( x_{k+1},\dots ,x_{n}\right)  \label{1.7}
\end{equation}%
and%
\begin{equation}
\Gamma ^{\frac{1}{2}}\left( x_{1}+y_{1},x_{2},\dots ,x_{n}\right) \leq
\Gamma ^{\frac{1}{2}}\left( x_{1},x_{2},\dots ,x_{n}\right) +\Gamma ^{\frac{1%
}{2}}\left( y_{1},x_{2},\dots ,x_{n}\right) .  \label{1.8}
\end{equation}

The main aim of this paper is to point out some upper bounds for the
distance $d\left( x,M\right) $ in terms of the linearly independent vectors $%
\left\{ x_{1},\dots ,x_{n}\right\} $ that span $M$ and $x\notin M^{\perp },$
where $M^{\perp }$ is the orthogonal complement of $M$ in the inner product
space $\left( H;\left\langle \cdot ,\cdot \right\rangle \right) $.

As a by-product of this endeavour, some refinements of the generalisations
for Bessel's inequality due to several authors including: Boas, Bellman and
Bombieri are obtained. Refinements for the well known Hadamard's inequality
for Gram determinants are also derived.

\section{Upper Bounds for $d\left( x,M\right) $}

The following result may be stated.

\begin{theorem}
\label{t2.1}Let $\left\{ x_{1},\dots ,x_{n}\right\} $ be a linearly
independent system of vectors in $H$ and $M:=span\left\{ x_{1},\dots
,x_{n}\right\} .$ If $x\notin M^{\perp },$ then%
\begin{equation}
d^{2}\left( x,M\right) <\frac{\left\Vert x\right\Vert
^{2}\sum_{i=1}^{n}\left\Vert x_{i}\right\Vert ^{2}-\sum_{i=1}^{n}\left\vert
\left\langle x,x_{i}\right\rangle \right\vert ^{2}}{\sum_{i=1}^{n}\left\Vert
x_{i}\right\Vert ^{2}}  \label{2.1}
\end{equation}%
or, equivalently,%
\begin{equation}
\Gamma \left( x_{1},\dots ,x_{n},x\right) <\frac{\left\Vert x\right\Vert
^{2}\sum_{i=1}^{n}\left\Vert x_{i}\right\Vert ^{2}-\sum_{i=1}^{n}\left\vert
\left\langle x,x_{i}\right\rangle \right\vert ^{2}}{\sum_{i=1}^{n}\left\Vert
x_{i}\right\Vert ^{2}}\cdot \Gamma \left( x_{1},\dots ,x_{n}\right) .
\label{2.2}
\end{equation}
\end{theorem}

\begin{proof}
If we use the Cauchy-Bunyakovsky-Schwarz type inequality%
\begin{equation}
\left\Vert \sum_{i=1}^{n}\alpha _{i}y_{i}\right\Vert ^{2}\leq
\sum_{i=1}^{n}\left\vert \alpha _{i}\right\vert ^{2}\sum_{i=1}^{n}\left\Vert
y_{i}\right\Vert ^{2},  \label{2.3}
\end{equation}%
that can be easily deduced from the obvious identity%
\begin{equation}
\sum_{i=1}^{n}\left\vert \alpha _{i}\right\vert ^{2}\sum_{i=1}^{n}\left\Vert
y_{i}\right\Vert ^{2}-\left\Vert \sum_{i=1}^{n}\alpha _{i}y_{i}\right\Vert
^{2}=\frac{1}{2}\sum_{i,j=1}^{n}\left\Vert \overline{\alpha _{i}}x_{j}-%
\overline{\alpha _{j}}x_{i}\right\Vert ^{2},  \label{2.4}
\end{equation}%
we can state that%
\begin{equation}
\left\Vert \sum_{i=1}^{n}\left\langle x,x_{i}\right\rangle x_{i}\right\Vert
^{2}\leq \sum_{i=1}^{n}\left\vert \left\langle x,x_{i}\right\rangle
\right\vert ^{2}\sum_{i=1}^{n}\left\Vert x_{i}\right\Vert ^{2}.  \label{2.5}
\end{equation}%
Note that the equality case holds in (\ref{2.5}) if and only if, by (\ref%
{2.4}), 
\begin{equation}
\overline{\left\langle x,x_{i}\right\rangle }x_{j}=\overline{\left\langle
x,x_{i}\right\rangle }x_{i}  \label{2.6}
\end{equation}%
for each $i,j\in \left\{ 1,\dots ,n\right\} .$

Utilising the expression (\ref{1.3}) of the distance $d\left( x,M\right) $,
we have%
\begin{equation}
d^{2}\left( x,M\right) =\left\Vert x\right\Vert ^{2}-\frac{%
\sum_{i=1}^{n}\left\vert \left\langle x,x_{i}\right\rangle \right\vert
^{2}\sum_{i=1}^{n}\left\Vert x_{i}\right\Vert ^{2}}{\left\Vert
\sum_{i=1}^{n}\left\langle x,x_{i}\right\rangle x_{i}\right\Vert ^{2}}\cdot 
\frac{\sum_{i=1}^{n}\left\vert \left\langle x,x_{i}\right\rangle \right\vert
^{2}}{\sum_{i=1}^{n}\left\Vert x_{i}\right\Vert ^{2}}.  \label{2.7}
\end{equation}%
Since $\left\{ x_{1},\dots ,x_{n}\right\} $ are linearly independent, hence (%
\ref{2.6}) cannot be achieved and then we have strict inequality in (\ref%
{2.5}).

Finally, on using (\ref{2.5}) and (\ref{2.7}) we get the desired result (\ref%
{2.1}).
\end{proof}

\begin{remark}
\label{r2.2}It is known that (see (\ref{1.4})) if not all $\left\{
x_{1},\dots ,x_{n}\right\} $ are orthogonal on each other, then the
following result which is well known in the literature as Hadamard's
inequality holds:%
\begin{equation}
\Gamma \left( x_{1},\dots ,x_{n}\right) <\left\Vert x_{1}\right\Vert
^{2}\left\Vert x_{2}\right\Vert ^{2}\cdots \left\Vert x_{n}\right\Vert ^{2}.
\label{2.8}
\end{equation}%
Utilising the inequality (\ref{2.2}), we may write successively:%
\begin{align*}
\Gamma \left( x_{1},x_{2}\right) & \leq \frac{\left\Vert x_{1}\right\Vert
^{2}\left\Vert x_{2}\right\Vert ^{2}-\left\vert \left\langle
x_{2},x_{1}\right\rangle \right\vert ^{2}}{\left\Vert x_{1}\right\Vert ^{2}}%
\left\Vert x_{1}\right\Vert ^{2}\leq \left\Vert x_{1}\right\Vert
^{2}\left\Vert x_{2}\right\Vert ^{2}, \\
\Gamma \left( x_{1},x_{2},x_{3}\right) & <\frac{\left\Vert x_{3}\right\Vert
^{2}\sum_{i=1}^{2}\left\Vert x_{i}\right\Vert ^{2}-\sum_{i=1}^{2}\left\vert
\left\langle x_{3},x_{i}\right\rangle \right\vert ^{2}}{\sum_{i=1}^{2}\left%
\Vert x_{i}\right\Vert ^{2}}\Gamma \left( x_{1},x_{2}\right) \\
& \leq \left\Vert x_{3}\right\Vert ^{2}\Gamma \left( x_{1},x_{2}\right) \\
& \cdots \cdots \cdots \cdots \cdots \cdots \cdots \cdots \cdots \cdots
\cdots \cdots \cdots \cdots \\
\Gamma \left( x_{1},\dots ,x_{n-1},x_{n}\right) & <\frac{\left\Vert
x_{n}\right\Vert ^{2}\sum_{i=1}^{n-1}\left\Vert x_{i}\right\Vert
^{2}-\sum_{i=1}^{n-1}\left\vert \left\langle x_{n},x_{i}\right\rangle
\right\vert ^{2}}{\sum_{i=1}^{n-1}\left\Vert x_{i}\right\Vert ^{2}}\Gamma
\left( x_{1},\dots ,x_{n-1}\right) \\
& \leq \left\Vert x_{n}\right\Vert ^{2}\Gamma \left( x_{1},\dots
,x_{n-1}\right) .
\end{align*}%
Multiplying the above inequalities, we deduce%
\begin{align}
\Gamma \left( x_{1},\dots ,x_{n-1},x_{n}\right) & <\left\Vert
x_{1}\right\Vert ^{2}\prod_{k=2}^{n}\left( \left\Vert x_{k}\right\Vert ^{2}-%
\frac{1}{\sum_{i=1}^{k-1}\left\Vert x_{i}\right\Vert ^{2}}%
\sum_{i=1}^{k-1}\left\vert \left\langle x_{k},x_{i}\right\rangle \right\vert
^{2}\right)  \label{2.9} \\
& \leq \prod_{j=1}^{n}\left\Vert x_{j}\right\Vert ^{2},  \notag
\end{align}%
valid for a system of $n\geq 2$ linearly independent vectors which are not
orthogonal on each other.
\end{remark}

In \cite{DRA1}, the author has obtained the following inequality.

\begin{lemma}
\label{l2.3}Let $z_{1},\dots ,z_{n}\in H$ and $\alpha _{1},\dots ,\alpha
_{n}\in \mathbb{K}$. Then one has the inequalities:%
\begin{multline}
\left\Vert \sum_{i=1}^{n}\alpha _{i}z_{i}\right\Vert ^{2}\leq \left\{ 
\begin{array}{l}
\max\limits_{1\leq i\leq n}\left\vert \alpha _{i}\right\vert
^{2}\sum\limits_{i=1}^{n}\left\Vert z_{i}\right\Vert ^{2}; \\ 
\\ 
\left( \sum\limits_{i=1}^{n}\left\vert \alpha _{i}\right\vert ^{2\alpha
}\right) ^{\frac{1}{\alpha }}\left( \sum\limits_{i=1}^{n}\left\Vert
z_{i}\right\Vert ^{2\beta }\right) ^{\frac{1}{p}} \\ 
\hfill \text{where \ }\alpha >1,\ \frac{1}{\alpha }+\frac{1}{\beta }=1; \\ 
\\ 
\sum\limits_{i=1}^{n}\left\vert \alpha _{i}\right\vert
^{2}\max\limits_{1\leq i\leq n}\left\Vert z_{i}\right\Vert ^{2};%
\end{array}%
\right.  \label{2.10} \\
+\left\{ 
\begin{array}{l}
\max\limits_{1\leq i\neq j\leq n}\left\{ \left\vert \alpha _{i}\alpha
_{j}\right\vert \right\} \sum\limits_{1\leq i\neq j\leq n}\left\vert
\left\langle z_{i},z_{j}\right\rangle \right\vert ; \\ 
\\ 
\left[ \left( \sum\limits_{i=1}^{n}\left\vert \alpha _{i}\right\vert
^{\gamma }\right) ^{2}-\sum\limits_{i=1}^{n}\left\vert \alpha
_{i}\right\vert ^{2\gamma }\right] ^{\frac{1}{\gamma }}\left(
\sum\limits_{1\leq i\neq j\leq n}\left\vert \left\langle
z_{i},z_{j}\right\rangle \right\vert ^{\delta }\right) ^{\frac{1}{\delta }}
\\ 
\hfill \text{where \ }\gamma >1,\ \frac{1}{\gamma }+\frac{1}{\delta }=1; \\ 
\\ 
\left[ \left( \sum\limits_{i=1}^{n}\left\vert \alpha _{i}\right\vert \right)
^{2}-\sum\limits_{i=1}^{n}\left\vert \alpha _{i}\right\vert ^{2}\right]
\max\limits_{1\leq i\neq j\leq n}\left\vert \left\langle
z_{i},z_{j}\right\rangle \right\vert ;%
\end{array}%
\right.
\end{multline}%
where any term in the first branch can be combined with each term from the
second branch giving 9 possible combinations.
\end{lemma}

Out of these, we select the following ones that are of relevance for further
consideration%
\begin{align}
& \left\Vert \sum_{i=1}^{n}\alpha _{i}z_{i}\right\Vert ^{2}  \label{2.11} \\
& \leq \max\limits_{1\leq i\leq n}\left\Vert z_{i}\right\Vert
^{2}\sum\limits_{i=1}^{n}\left\vert \alpha _{i}\right\vert
^{2}+\max\limits_{1\leq i<j\leq n}\left\vert \left\langle
z_{i},z_{j}\right\rangle \right\vert \left[ \left(
\sum\limits_{i=1}^{n}\left\vert \alpha _{i}\right\vert \right)
^{2}-\sum\limits_{i=1}^{n}\left\vert \alpha _{i}\right\vert ^{2}\right] 
\notag \\
& \leq \sum\limits_{i=1}^{n}\left\vert \alpha _{i}\right\vert ^{2}\left(
\max\limits_{1\leq i\leq n}\left\Vert z_{i}\right\Vert ^{2}+\left(
n-1\right) \max\limits_{1\leq i<j\leq n}\left\vert \left\langle
z_{i},z_{j}\right\rangle \right\vert \right)  \notag
\end{align}%
and%
\begin{align}
& \left\Vert \sum_{i=1}^{n}\alpha _{i}z_{i}\right\Vert ^{2}  \label{2.12} \\
& \leq \max\limits_{1\leq i\leq n}\left\Vert z_{i}\right\Vert
^{2}\sum\limits_{i=1}^{n}\left\vert \alpha _{i}\right\vert ^{2}+\left[
\left( \sum\limits_{i=1}^{n}\left\vert \alpha _{i}\right\vert ^{2}\right)
^{2}-\sum\limits_{i=1}^{n}\left\vert \alpha _{i}\right\vert ^{4}\right]
^{1/2}  \notag \\
& \qquad \qquad \qquad \qquad \qquad \qquad \qquad \times \left(
\sum\limits_{1\leq i\neq j\leq n}\left\vert \left\langle
z_{i},z_{j}\right\rangle \right\vert ^{2}\right) ^{\frac{1}{2}}  \notag \\
& \leq \sum\limits_{i=1}^{n}\left\vert \alpha _{i}\right\vert ^{2}\left[
\max\limits_{1\leq i\leq n}\left\Vert z_{i}\right\Vert ^{2}+\left(
\sum\limits_{1\leq i\neq j\leq n}\left\vert \left\langle
z_{i},z_{j}\right\rangle \right\vert ^{2}\right) ^{\frac{1}{2}}\right] . 
\notag
\end{align}%
Note that the last inequality in (\ref{2.11}) follows by the fact that%
\begin{equation*}
\left( \sum\limits_{i=1}^{n}\left\vert \alpha _{i}\right\vert \right)
^{2}\leq n\sum\limits_{i=1}^{n}\left\vert \alpha _{i}\right\vert ^{2},
\end{equation*}%
while the last inequality in (\ref{2.12}) is obvious.

Utilising the above inequalities (\ref{2.11}) and (\ref{2.12}) which provide
alternatives to the Cauchy-Bunyakovsky-Schwarz inequality (\ref{2.3}), we
can state the following results.

\begin{theorem}
\label{t2.4}Let $\left\{ x_{1},\dots ,x_{n}\right\} ,$ $M$ and $x$ be as in
Theorem \ref{t2.1}. Then%
\begin{equation}
d^{2}\left( x,M\right) \leq \frac{\left\Vert x\right\Vert ^{2}\left[
\max\limits_{1\leq i\leq n}\left\Vert x_{i}\right\Vert ^{2}+\left(
\sum\limits_{1\leq i\neq j\leq n}\left\vert \left\langle
x_{i},x_{j}\right\rangle \right\vert ^{2}\right) ^{\frac{1}{2}}\right]
-\sum\limits_{i=1}^{n}\left\vert \left\langle x,x_{i}\right\rangle
\right\vert ^{2}}{\max\limits_{1\leq i\leq n}\left\Vert x_{i}\right\Vert
^{2}+\left( \sum\limits_{1\leq i\neq j\leq n}\left\vert \left\langle
x_{i},x_{j}\right\rangle \right\vert ^{2}\right) ^{\frac{1}{2}}}
\label{2.13}
\end{equation}%
or, equivalently,%
\begin{multline}
\Gamma \left( x_{1},\dots ,x_{n},x\right)  \label{2.14} \\
\leq \frac{\left\Vert x\right\Vert ^{2}\left[ \max\limits_{1\leq i\leq
n}\left\Vert x_{i}\right\Vert ^{2}+\left( \sum\limits_{1\leq i\neq j\leq
n}\left\vert \left\langle x_{i},x_{j}\right\rangle \right\vert ^{2}\right) ^{%
\frac{1}{2}}\right] -\sum\limits_{i=1}^{n}\left\vert \left\langle
x,x_{i}\right\rangle \right\vert ^{2}}{\max\limits_{1\leq i\leq n}\left\Vert
x_{i}\right\Vert ^{2}+\left( \sum\limits_{1\leq i\neq j\leq n}\left\vert
\left\langle x_{i},x_{j}\right\rangle \right\vert ^{2}\right) ^{\frac{1}{2}}}
\\
\times \Gamma \left( x_{1},\dots ,x_{n}\right)
\end{multline}
\end{theorem}

\begin{proof}
Utilising the inequality (\ref{2.12}) for $\alpha _{i}=\left\langle
x,x_{i}\right\rangle $ and $z_{i}=x_{i},$ $i\in \left\{ 1,\dots ,n\right\} ,$
we can write:%
\begin{equation}
\left\Vert \sum_{i=1}^{n}\left\langle x,x_{i}\right\rangle x_{i}\right\Vert
^{2}\leq \sum_{i=1}^{n}\left\vert \left\langle x,x_{i}\right\rangle
\right\vert ^{2}\left[ \max\limits_{1\leq i\leq n}\left\Vert
x_{i}\right\Vert ^{2}+\left( \sum_{1\leq i\neq j\leq n}\left\vert
\left\langle x_{i},x_{j}\right\rangle \right\vert ^{2}\right) ^{\frac{1}{2}}%
\right]  \label{2.15}
\end{equation}%
for any $x\in H.$

Now, since, by the representation formula (\ref{1.3})%
\begin{equation}
d^{2}\left( x,M\right) =\left\Vert x\right\Vert ^{2}-\frac{%
\sum_{i=1}^{n}\left\vert \left\langle x,x_{i}\right\rangle \right\vert ^{2}}{%
\left\Vert \sum_{i=1}^{n}\left\langle x,x_{i}\right\rangle x_{i}\right\Vert
^{2}}\cdot \sum_{i=1}^{n}\left\vert \left\langle x,x_{i}\right\rangle
\right\vert ^{2},  \label{2.16}
\end{equation}%
for $x\notin M^{\perp },$ hence, by (\ref{2.15}) and (\ref{2.16}) we deduce
the desired result (\ref{2.13}).
\end{proof}

\begin{remark}
\label{r2.5}In 1941, R.P. Boas \cite{BO} and in 1944, R. Bellman \cite{BE},
independent of each other, proved the following generalisation of Bessel's
inequality:%
\begin{equation}
\sum_{i=1}^{n}\left\vert \left\langle y,y_{i}\right\rangle \right\vert
^{2}\leq \left\Vert y\right\Vert ^{2}\left[ \max\limits_{1\leq i\leq
n}\left\Vert y_{i}\right\Vert ^{2}+\left( \sum_{1\leq i\neq j\leq
n}\left\vert \left\langle y_{i},y_{j}\right\rangle \right\vert ^{2}\right) ^{%
\frac{1}{2}}\right] ,  \label{2.17}
\end{equation}%
provided $y$ and $y_{i}$ $\left( i\in \left\{ 1,\dots ,n\right\} \right) $
are arbitrary vectors in the inner product space $\left( H;\left\langle
\cdot ,\cdot \right\rangle \right) .$ If $\left\{ y_{i}\right\} _{i\in
\left\{ 1,\dots ,n\right\} }$ are orthonormal, then (\ref{2.17}) reduces to
Bessel's inequality.

In this respect, one may see (\ref{2.13}) as a refinement of the
Boas-Bellman result (\ref{2.17}).
\end{remark}

\begin{remark}
\label{r2.6}On making use of a similar argument to that utilised in Remark %
\ref{r2.2}, one can obtain the following refinement of the Hadamard
inequality:%
\begin{align}
& \Gamma \left( x_{1},\dots ,x_{n}\right)  \label{2.18} \\
& \leq \left\Vert x_{1}\right\Vert ^{2}\prod_{k=2}^{n}\left( \left\Vert
x_{k}\right\Vert ^{2}-\frac{\sum\limits_{i=1}^{k-1}\left\vert \left\langle
x_{k},x_{i}\right\rangle \right\vert ^{2}}{\max\limits_{1\leq i\leq
k-1}\left\Vert x_{i}\right\Vert ^{2}+\left( \sum\limits_{1\leq i\neq j\leq
k-1}\left\vert \left\langle x_{i},x_{j}\right\rangle \right\vert ^{2}\right)
^{\frac{1}{2}}}\right)  \notag \\
& \leq \prod_{j=1}^{n}\left\Vert x_{j}\right\Vert ^{2}.  \notag
\end{align}
\end{remark}

Further on, if we choose $\alpha _{i}=\left\langle x,x_{i}\right\rangle ,$ $%
z_{i}=x_{i},$ $i\in \left\{ 1,\dots ,n\right\} $ in (\ref{2.11}), then we
may state the inequality%
\begin{equation}
\left\Vert \sum_{i=1}^{n}\left\langle x,x_{i}\right\rangle x_{i}\right\Vert
^{2}\leq \sum_{i=1}^{n}\left\vert \left\langle x,x_{i}\right\rangle
\right\vert ^{2}\left( \max\limits_{1\leq i\leq n}\left\Vert
x_{i}\right\Vert ^{2}+\left( n-1\right) \max_{1\leq i\neq j\leq n}\left\vert
\left\langle x_{i},x_{j}\right\rangle \right\vert \right) .  \label{2.19}
\end{equation}%
Utilising (\ref{2.19}) and (\ref{2.16}) we may state the following result as
well:

\begin{theorem}
\label{t2.7}Let $\left\{ x_{1},\dots ,x_{n}\right\} ,$ $M$ and $x$ be as in
Theorem \ref{t2.1}. Then%
\begin{multline}
d^{2}\left( x,M\right)  \label{2.20} \\
\leq \frac{\left\Vert x\right\Vert ^{2}\left[ \max\limits_{1\leq i\leq
n}\left\Vert x_{i}\right\Vert ^{2}+\left( n-1\right) \max\limits_{1\leq
i\neq j\leq n}\left\vert \left\langle x_{i},x_{j}\right\rangle \right\vert %
\right] -\sum_{i=1}^{n}\left\vert \left\langle x,x_{i}\right\rangle
\right\vert ^{2}}{\max\limits_{1\leq i\leq n}\left\Vert x_{i}\right\Vert
^{2}+\left( n-1\right) \max\limits_{1\leq i\neq j\leq n}\left\vert
\left\langle x_{i},x_{j}\right\rangle \right\vert }
\end{multline}%
or, equivalently,%
\begin{multline}
\Gamma \left( x_{1},\dots ,x_{n},x\right)  \label{2.21} \\
\leq \frac{\left\Vert x\right\Vert ^{2}\left[ \max\limits_{1\leq i\leq
n}\left\Vert x_{i}\right\Vert ^{2}+\left( n-1\right) \max\limits_{1\leq
i\neq j\leq n}\left\vert \left\langle x_{i},x_{j}\right\rangle \right\vert %
\right] -\sum_{i=1}^{n}\left\vert \left\langle x,x_{i}\right\rangle
\right\vert ^{2}}{\max\limits_{1\leq i\leq n}\left\Vert x_{i}\right\Vert
^{2}+\left( n-1\right) \max\limits_{1\leq i\neq j\leq n}\left\vert
\left\langle x_{i},x_{j}\right\rangle \right\vert } \\
\times \Gamma \left( x_{1},\dots ,x_{n}\right)
\end{multline}
\end{theorem}

\begin{remark}
\label{r2.8}The above result (\ref{2.20}) provides a refinement for the
following generalisation of Bessel's inequality:%
\begin{equation}
\sum_{i=1}^{n}\left\vert \left\langle x,x_{i}\right\rangle \right\vert
^{2}\leq \left\Vert x\right\Vert ^{2}\left[ \max\limits_{1\leq i\leq
n}\left\Vert x_{i}\right\Vert ^{2}+\left( n-1\right) \max\limits_{1\leq
i\neq j\leq n}\left\vert \left\langle x_{i},x_{j}\right\rangle \right\vert %
\right] ,  \label{2.22}
\end{equation}%
obtained by the author in \cite{DRA1}.

One can also provide the corresponding refinement of Hadamard's inequality (%
\ref{1.4}) on using (\ref{2.21}), i.e., 
\begin{align}
& \Gamma \left( x_{1},\dots ,x_{n}\right)  \label{2.23} \\
& \leq \left\Vert x_{1}\right\Vert ^{2}\prod_{k=2}^{n}\left( \left\Vert
x_{k}\right\Vert ^{2}-\frac{\sum\limits_{i=1}^{k-1}\left\vert \left\langle
x_{k},x_{i}\right\rangle \right\vert ^{2}}{\max\limits_{1\leq i\leq
k-1}\left\Vert x_{i}\right\Vert ^{2}+\left( k-2\right) \max\limits_{1\leq
i\neq j\leq k-1}\left\vert \left\langle x_{i},x_{j}\right\rangle \right\vert 
}\right)  \notag \\
& \leq \prod_{j=1}^{n}\left\Vert x_{j}\right\Vert ^{2}.  \notag
\end{align}
\end{remark}

\section{Other Upper Bounds for $d\left( x,M\right) $}

In \cite[p. 140]{DRA2} the author obtained the following inequality that is
similar to the Cauchy-Bunyakovsky-Schwarz result.

\begin{lemma}
\label{l3.1}Let $z_{1},\dots ,z_{n}\in H$ and $\alpha _{1},\dots ,\alpha
_{n}\in \mathbb{K}$. Then one has the inequalities:%
\begin{align}
\left\Vert \sum_{i=1}^{n}\alpha _{i}z_{i}\right\Vert ^{2}& \leq
\sum\limits_{i=1}^{n}\left\vert \alpha _{i}\right\vert
^{2}\sum\limits_{j=1}^{n}\left\vert \left\langle z_{i},z_{j}\right\rangle
\right\vert  \label{3.1} \\
& \leq \left\{ 
\begin{array}{l}
\sum\limits_{i=1}^{n}\left\vert \alpha _{i}\right\vert
^{2}\max\limits_{1\leq i\leq n}\left[ \sum\limits_{j=1}^{n}\left\vert
\left\langle z_{i},z_{j}\right\rangle \right\vert \right] ; \\ 
\\ 
\left( \sum\limits_{i=1}^{n}\left\vert \alpha _{i}\right\vert ^{2p}\right) ^{%
\frac{1}{p}}\left( \sum\limits_{i=1}^{n}\left(
\sum\limits_{j=1}^{n}\left\vert \left\langle z_{i},z_{j}\right\rangle
\right\vert \right) ^{q}\right) ^{\frac{1}{q}} \\ 
\hfill \text{where \ }p>1,\ \frac{1}{p}+\frac{1}{q}=1; \\ 
\\ 
\max\limits_{1\leq i\leq n}\left\vert \alpha _{i}\right\vert
^{2}\sum\limits_{i,j=1}^{n}\left\vert \left\langle z_{i},z_{j}\right\rangle
\right\vert .%
\end{array}%
\right.  \notag
\end{align}
\end{lemma}

We can state and prove now another upper bound for the distance $d\left(
x,M\right) $ as follows.

\begin{theorem}
\label{t3.2}Let $\left\{ x_{1},\dots ,x_{n}\right\} ,$ $M$ and $x$ be as in
Theorem \ref{t2.1}. Then%
\begin{equation}
d^{2}\left( x,M\right) \leq \frac{\left\Vert x\right\Vert
^{2}\max\limits_{1\leq i\leq n}\left[ \sum\limits_{j=1}^{n}\left\vert
\left\langle x_{i},x_{j}\right\rangle \right\vert \right] -\sum%
\limits_{i=1}^{n}\left\vert \left\langle x,x_{i}\right\rangle \right\vert
^{2}}{\max\limits_{1\leq i\leq n}\left[ \sum\limits_{j=1}^{n}\left\vert
\left\langle x_{i},x_{j}\right\rangle \right\vert \right] }  \label{3.2}
\end{equation}%
or, equivalently,%
\begin{multline}
\Gamma \left( x_{1},\dots ,x_{n},x\right)  \label{3.3} \\
\leq \frac{\left\Vert x\right\Vert ^{2}\max\limits_{1\leq i\leq n}\left[
\sum\limits_{j=1}^{n}\left\vert \left\langle x_{i},x_{j}\right\rangle
\right\vert \right] -\sum\limits_{i=1}^{n}\left\vert \left\langle
x,x_{i}\right\rangle \right\vert ^{2}}{\max\limits_{1\leq i\leq n}\left[
\sum\limits_{j=1}^{n}\left\vert \left\langle x_{i},x_{j}\right\rangle
\right\vert \right] }\cdot \Gamma \left( x_{1},\dots ,x_{n}\right) .
\end{multline}
\end{theorem}

\begin{proof}
Utilising the first branch in (\ref{3.1}) we may state that%
\begin{equation}
\left\Vert \sum_{i=1}^{n}\left\langle x,x_{i}\right\rangle x_{i}\right\Vert
^{2}\leq \sum_{i=1}^{n}\left\vert \left\langle x,x_{i}\right\rangle
\right\vert ^{2}\max\limits_{1\leq i\leq n}\left[ \sum\limits_{j=1}^{n}\left%
\vert \left\langle x_{i},x_{j}\right\rangle \right\vert \right]  \label{3.3a}
\end{equation}%
for any $x\in H.$

Now, since, by the representation formula (\ref{1.3}) we have%
\begin{equation}
d^{2}\left( x,M\right) =\left\Vert x\right\Vert ^{2}-\frac{%
\sum_{i=1}^{n}\left\vert \left\langle x,x_{i}\right\rangle \right\vert ^{2}}{%
\left\Vert \sum_{i=1}^{n}\left\langle x,x_{i}\right\rangle x_{i}\right\Vert
^{2}}\cdot \sum_{i=1}^{n}\left\vert \left\langle x,x_{i}\right\rangle
\right\vert ^{2},  \label{3.4}
\end{equation}%
for $x\notin M^{\perp },$ hence, by (\ref{3.3a}) and (\ref{3.4}) we deduce
the desired result (\ref{3.2}).
\end{proof}

\begin{remark}
\label{r3.3}In 1971, E. Bombieri \cite{BOM} proved the following
generalisation of Bessel's inequality, however not stated in the general
form for inner products. The general version can be found for instance in 
\cite[p. 394]{MPF}. It reads as follows: if $y,y_{1},\dots ,y_{n}$ are
vectors in the inner product space $\left( H;\left\langle \cdot ,\cdot
\right\rangle \right) ,$ then%
\begin{equation}
\sum_{i=1}^{n}\left\vert \left\langle y,y_{i}\right\rangle \right\vert
^{2}\leq \left\Vert y\right\Vert ^{2}\max\limits_{1\leq i\leq n}\left\{
\sum\limits_{j=1}^{n}\left\vert \left\langle y_{i},y_{j}\right\rangle
\right\vert \right\} .  \label{3.5}
\end{equation}%
Obviously, when $\left\{ y_{1},\dots ,y_{n}\right\} $ are orthonormal, the
inequality (\ref{3.5}) produces Bessel's inequality.

In this respect, we may regard our result (\ref{3.2}) as a refinement of the
Bombieri inequality (\ref{3.5}).
\end{remark}

\begin{remark}
\label{r3.4}On making use of a similar argument to that in Remark \ref{r2.2}%
, we obtain the following refinement for the Hadamard inequality:%
\begin{align}
\Gamma \left( x_{1},\dots ,x_{n}\right) & \leq \left\Vert x_{1}\right\Vert
^{2}\prod_{k=2}^{n}\left[ \left\Vert x_{k}\right\Vert ^{2}-\frac{%
\sum\limits_{i=1}^{k-1}\left\vert \left\langle x_{k},x_{i}\right\rangle
\right\vert ^{2}}{\max\limits_{1\leq i\leq k-1}\left[ \sum%
\limits_{j=1}^{k-1}\left\vert \left\langle x_{i},x_{j}\right\rangle
\right\vert \right] }\right]  \label{3.6} \\
& \leq \prod_{j=1}^{n}\left\Vert x_{j}\right\Vert ^{2}.  \notag
\end{align}
\end{remark}

Another different Cauchy-Bunyakovsky-Schwarz type inequality is incorporated
in the following lemma \cite{DRA3}.

\begin{lemma}
\label{l3.4}Let $z_{1},\dots ,z_{n}\in H$ and $\alpha _{1},\dots ,\alpha
_{n}\in \mathbb{K}$. Then 
\begin{equation}
\left\Vert \sum_{i=1}^{n}\alpha _{i}z_{i}\right\Vert ^{2}\leq \left(
\sum\limits_{i=1}^{n}\left\vert \alpha _{i}\right\vert ^{p}\right) ^{\frac{2%
}{p}}\left( \sum\limits_{i,j=1}^{n}\left\vert \left\langle
z_{i},z_{j}\right\rangle \right\vert ^{q}\right) ^{\frac{1}{q}}  \label{3.7}
\end{equation}%
for $p>1,$ $\frac{1}{p}+\frac{1}{q}=1.$

If in (\ref{3.7}) we choose $p=q=2,$ then we get%
\begin{equation}
\left\Vert \sum_{i=1}^{n}\alpha _{i}z_{i}\right\Vert ^{2}\leq
\sum\limits_{i=1}^{n}\left\vert \alpha _{i}\right\vert ^{2}\left(
\sum\limits_{i,j=1}^{n}\left\vert \left\langle z_{i},z_{j}\right\rangle
\right\vert ^{2}\right) ^{\frac{1}{2}}.  \label{3.8}
\end{equation}
\end{lemma}

Based on (\ref{3.8}), we can state the following result that provides yet
another upper bound for the distance $d\left( x,M\right) .$

\begin{theorem}
\label{t3.5}Let $\left\{ x_{1},\dots ,x_{n}\right\} ,$ $M$ and $x$ be as in
Theorem \ref{t2.1}. Then%
\begin{equation}
d^{2}\left( x,M\right) \leq \frac{\left\Vert x\right\Vert ^{2}\left(
\sum\limits_{i,j=1}^{n}\left\vert \left\langle x_{i},x_{j}\right\rangle
\right\vert ^{2}\right) ^{\frac{1}{2}}-\sum\limits_{i=1}^{n}\left\vert
\left\langle x,x_{i}\right\rangle \right\vert ^{2}}{\left(
\sum\limits_{i,j=1}^{n}\left\vert \left\langle x_{i},x_{j}\right\rangle
\right\vert ^{2}\right) ^{\frac{1}{2}}}  \label{3.9}
\end{equation}%
or, equivalently,%
\begin{multline}
\Gamma \left( x_{1},\dots ,x_{n},x\right)  \label{3.10} \\
\leq \frac{\left\Vert x\right\Vert ^{2}\left(
\sum\limits_{i,j=1}^{n}\left\vert \left\langle x_{i},x_{j}\right\rangle
\right\vert ^{2}\right) ^{\frac{1}{2}}-\sum\limits_{i=1}^{n}\left\vert
\left\langle x,x_{i}\right\rangle \right\vert ^{2}}{\left(
\sum\limits_{i,j=1}^{n}\left\vert \left\langle x_{i},x_{j}\right\rangle
\right\vert ^{2}\right) ^{\frac{1}{2}}}\cdot \Gamma \left( x_{1},\dots
,x_{n}\right) .
\end{multline}
\end{theorem}

Similar comments apply related to Hadamard's inequality. We omit the details.

\section{Some Conditional Bounds}

In the recent paper \cite{DRA4}, the author has established the following
reverse of the Bessel inequality.

Let $\left( H;\left\langle \cdot ,\cdot \right\rangle \right) $ be an inner
product space over the real or complex number field $\mathbb{K}$, $\left\{
e_{i}\right\} _{i\in I}$ a finite family of orthonormal vectors in $H,$ $%
\varphi _{i},\phi _{i}\in \mathbb{K}$, $i\in I$ and $x\in H.$ If%
\begin{equation}
\func{Re}\left\langle \sum_{i\in I}\phi _{i}e_{i}-x,x-\sum_{i\in I}\varphi
_{i}e_{i}\right\rangle \geq 0  \label{4.1}
\end{equation}%
or, equivalently,%
\begin{equation}
\left\Vert x-\sum_{i\in I}\frac{\varphi _{i}+\phi _{i}}{2}e_{i}\right\Vert
\leq \frac{1}{2}\left( \sum_{i\in I}\left\vert \phi _{i}-\varphi
_{i}\right\vert ^{2}\right) ^{\frac{1}{2}},  \label{4.2}
\end{equation}%
then%
\begin{equation}
\left( 0\leq \right) \left\Vert x\right\Vert ^{2}-\sum_{i\in I}\left\vert
\left\langle x,e_{i}\right\rangle \right\vert ^{2}\leq \frac{1}{4}\sum_{i\in
I}\left\vert \phi _{i}-\varphi _{i}\right\vert ^{2}.  \label{4.3}
\end{equation}%
The constant $\frac{1}{4}$ is best possible in the sense that it cannot be
replaced by a smaller constant.

\begin{theorem}
\label{t4.1}Let $\left\{ x_{1},\dots x_{n}\right\} $ be a linearly
independent system of vectors in $H$ and $M:=span\left\{ x_{1},\dots
x_{n}\right\} .$ If $\gamma _{i},$ $\Gamma _{i}\in \mathbb{K}$, $i\in
\left\{ 1,\dots ,n\right\} $ and $x\in H\backslash M^{\perp }$ is such that%
\begin{equation}
\func{Re}\left\langle \sum_{i=1}^{n}\Gamma
_{i}x_{i}-x,x-\sum_{i=1}^{n}\gamma _{i}x_{i}\right\rangle \geq 0,
\label{4.4}
\end{equation}%
then we have the bound%
\begin{equation}
d^{2}\left( x,M\right) \leq \frac{1}{4}\left\Vert \sum_{i=1}^{n}\left(
\Gamma _{i}-\gamma _{i}\right) x_{i}\right\Vert ^{2}  \label{4.5}
\end{equation}%
or, equivalently,%
\begin{equation}
\Gamma \left( x_{1},\dots ,x_{n},x\right) \leq \frac{1}{4}\left\Vert
\sum_{i=1}^{n}\left( \Gamma _{i}-\gamma _{i}\right) x_{i}\right\Vert
^{2}\Gamma \left( x_{1},\dots ,x_{n}\right) .  \label{4.6}
\end{equation}
\end{theorem}

\begin{proof}
It is easy to see that in an inner product space for any $x,z,Z\in H$ one has%
\begin{equation*}
\left\Vert x-\frac{z+Z}{2}\right\Vert ^{2}-\frac{1}{4}\left\Vert
Z-z\right\Vert ^{2}=\func{Re}\left\langle Z-x,x-z\right\rangle ,
\end{equation*}%
therefore, the condition (\ref{4.4}) is actually equivalent to%
\begin{equation}
\left\Vert x-\sum_{i=1}^{n}\frac{\Gamma _{i}+\gamma _{i}}{2}x_{i}\right\Vert
^{2}\leq \frac{1}{4}\left\Vert \sum_{i=1}^{n}\left( \Gamma _{i}-\gamma
_{i}\right) x_{i}\right\Vert ^{2}.  \label{4.7}
\end{equation}%
Now, obviously,%
\begin{equation}
d^{2}\left( x,M\right) =\inf_{y\in M}\left\Vert x-y\right\Vert ^{2}\leq
\left\Vert x-\sum_{i=1}^{n}\frac{\Gamma _{i}+\gamma _{i}}{2}x_{i}\right\Vert
^{2}  \label{4.8}
\end{equation}%
and thus, by (\ref{4.7}) and (\ref{4.8}) we deduce (\ref{4.5}).

The last inequality is obvious by the representation (\ref{1.2}).
\end{proof}

\begin{remark}
\label{r4.2}Utilising various Cauchy-Bunyakovsky-Schwarz type inequalities
we may obtain more convenient (although coarser) bounds for $d^{2}\left(
x,M\right) .$ For instance, if we use the inequality (\ref{2.11}) we can
state the inequality:%
\begin{equation*}
\left\Vert \sum_{i=1}^{n}\left( \Gamma _{i}-\gamma _{i}\right)
x_{i}\right\Vert ^{2}\leq \sum_{i=1}^{n}\left\vert \Gamma _{i}-\gamma
_{i}\right\vert ^{2}\left( \max\limits_{1\leq i\leq n}\left\Vert
x_{i}\right\Vert ^{2}+\left( n-1\right) \max\limits_{1\leq i<j\leq
n}\left\vert \left\langle x_{i},x_{j}\right\rangle \right\vert \right) ,
\end{equation*}%
giving the bound:%
\begin{equation}
d^{2}\left( x,M\right) \leq \frac{1}{4}\sum_{i=1}^{n}\left\vert \Gamma
_{i}-\gamma _{i}\right\vert ^{2}\left[ \max\limits_{1\leq i\leq n}\left\Vert
x_{i}\right\Vert ^{2}+\left( n-1\right) \max\limits_{1\leq i<j\leq
n}\left\vert \left\langle x_{i},x_{j}\right\rangle \right\vert \right] ,
\label{4.9}
\end{equation}%
provided (\ref{4.4}) holds true.

Obviously, if $\left\{ x_{1},\dots ,x_{n}\right\} $ is an orthonormal family
in $H,$ then from (\ref{4.9}) we deduce the reverse of Bessel's inequality
incorporated in (\ref{4.3}).

If we use the inequality (\ref{2.12}), then we can state the inequality%
\begin{equation*}
\left\Vert \sum_{i=1}^{n}\left( \Gamma _{i}-\gamma _{i}\right)
x_{i}\right\Vert ^{2}\leq \sum_{i=1}^{n}\left\vert \Gamma _{i}-\gamma
_{i}\right\vert ^{2}\left[ \max\limits_{1\leq i\leq n}\left\Vert
x_{i}\right\Vert ^{2}+\left( \sum\limits_{1\leq i\neq j\leq n}\left\vert
\left\langle x_{i},x_{j}\right\rangle \right\vert ^{2}\right) ^{\frac{1}{2}}%
\right] ,
\end{equation*}%
giving the bound%
\begin{equation}
d^{2}\left( x,M\right) \leq \frac{1}{4}\sum_{i=1}^{n}\left\vert \Gamma
_{i}-\gamma _{i}\right\vert ^{2}\left[ \max\limits_{1\leq i\leq n}\left\Vert
x_{i}\right\Vert ^{2}+\left( \sum\limits_{1\leq i\neq j\leq n}\left\vert
\left\langle x_{i},x_{j}\right\rangle \right\vert ^{2}\right) ^{\frac{1}{2}}%
\right] ,  \label{4.10}
\end{equation}%
provided (\ref{4.4}) holds true.

In this case, when one assumes that $\left\{ x_{1},\dots ,x_{n}\right\} $ is
an orthonormal family of vectors, then (\ref{4.10}) reduces to (\ref{4.3})
as well.

Finally, on utilising the first branch of the inequality (\ref{3.1}), we can
state that%
\begin{equation}
d^{2}\left( x,M\right) \leq \frac{1}{4}\sum_{i=1}^{n}\left\vert \Gamma
_{i}-\gamma _{i}\right\vert ^{2}\max\limits_{1\leq i\leq n}\left[
\sum_{j=1}^{n}\left\vert \left\langle x_{i},x_{j}\right\rangle \right\vert %
\right] ,  \label{4.11}
\end{equation}%
provided (\ref{4.4}) holds true.

This inequality is also a generalisation of (\ref{4.3}).
\end{remark}


\begin{thebibliography}{9}
\bibitem{BE} R. BELLMAN, Almost orthogonal series,\textit{\ Bull. Amer.
Math. Soc}., \textbf{50} (1944), 517-519.

\bibitem{BO} R.P. BOAS, A general moment problem,\textit{\ Amer. J. Math}., 
\textbf{63} (1941), 361-370.

\bibitem{BOM} E. BOMBIERI, A note on the large sieve,\textit{\ Acta Arith}., 
\textbf{18} (1971), 401-404.

\bibitem{DE} F. DEUTSCH, \textit{Best Approximation in Inner Product Spaces, 
}CMS Books in Mathematics, Springer-Verlag, New York, Berlin, Heidelberg,
2001.

\bibitem{DRA4} S.S. DRAGOMIR, A counterpart of Bessel's inequality in inner
product spaces and some Gr\"{u}ss type related results, \textit{RGMIA Res.
Rep. Coll., }\textbf{6} (2003), Supplement, Article 10. [ONLINE: \texttt{http://rgmia.vu.edu.au/v6(E).html}].

\bibitem{DRA3} S.S. DRAGOMIR, Some Bombieri type inequalities in inner
product spaces,\textit{\ J. Indones. Math. Soc}., \textbf{10}(2) (2004),
91-97.

\bibitem{DRA1} S.S. DRAGOMIR, On the Boas-Bellman inequality in inner
product spaces,\textit{\ Bull. Austral. Math. Soc}., \textbf{69}(2) (2004),
217-225.

\bibitem{DRA2} S.S. DRAGOMIR, \textit{Advances in Inequalities of the
Schwarz, Gr\"{u}ss and Bessel Type in Inner Product Spaces}, RGMIA
Monographs, Victoria University, 2004. [ONLINE: \texttt{%
http://rgmia.vu.edu.au/monographs/}].

\bibitem{MPF} D.S. MITRINOVI\'{C}, J.E. PE\v{C}ARI\'{C} and\ A.M. FINK, 
\textit{Classical and New Inequalities in Analysis, }Kluwer Academic,
Dordrecht, 1993.
\end{thebibliography}
\end{document}